\documentclass{amsart}

\usepackage{amsfonts}
\usepackage{latexsym}
\usepackage{amssymb}
\usepackage{amscd}

\newtheorem{thm}{Theorem}[section]

\newtheorem{prop}[thm]{Proposition}
\newtheorem{cor}[thm]{Corollary}

\newcommand{\A}{{\mathcal{A}}}
\newcommand{\B}{{\mathcal{B}}}

\def\qed{{$\Box$}\medskip}

  \def\s{\sigma}

\def\A{\mathcal{A}}

\def\C{\mathcal{C}}

\def\H{\mathrm{H}}

\def\M{\mathcal{M}}
\def\N{\mathcal{N}}

\def\U{\mathcal{U}}

\def\X{\mathcal{X}}
\def\Y{\mathcal{Y}}

\def\Ext{\mathrm{Ext}}

\def\slq2{SL_q (2)}
\def\s{\sigma}
\def\smod{\sigma}

\def\hu{\hat{u}}

\def\HH{\mathrm{HH}}
\def\Hdim{\mathrm{dim}}

\def\gldim{\mathrm{gldim}}

\begin{document}
\title{On the Hochschild homology of quantum $SL(N)$}
\author{Tom~Hadfield${}^1$, Ulrich~Kr\"{a}hmer${}^2$} 
\date{\today}
\maketitle

\centerline{${}^1$ School of Mathematical Sciences, Queen Mary, University of London}
\centerline{327 Mile End Road, London E1 4NS, England, t.hadfield@qmul.ac.uk}
\centerline{Supported  by an EPSRC postdoctoral fellowship.}
\centerline{}
\centerline{ ${}^2$ Instytut Matematyczny Polskiej Akademii Nauk}
\centerline{Ul. Sniadeckich 8, PL-00956 Warszawa, Poland, kraehmer@impan.gov.pl}
\centerline{Supported by an EU Marie Curie postdoctoral fellowship.}
\centerline{}
\centerline{MSC (2000): 16E40, 16W35}

\begin{abstract}
We show that the quantized coordinate ring $\A:=k_q[SL(N)]$ 
satisfies van den Bergh's analogue of Poincar\'e duality for Hochschild
(co)homology with dualizing bimodule being $\A_\sigma$, the 
$\A$-bimodule
which is $\A$ as $k$-vector space with right
multiplication twisted by the modular automorphism $ \sigma $ of the
Haar functional. This implies that 
$\H_{N^2-1}(\A,\A_\sigma) \cong k$,
generalizing our previous result for $k_q [ SL(2)]$.
\end{abstract}

\section{Introduction and statement of the result} 
\label{}
According to the 
Hochschild-Kostant-Rosenberg theorem 
\cite{aHKR}, the dimension of a regular affine variety $V$ 
over an algebraically closed field $k$ of characteristic 0 can be
expressed in terms of the Hochschild homology of its
coordinate ring $k[V]$ as 
\begin{equation}
\label{hoch_dim}
		\dim(V)=\mathrm{sup}\{n \ge 0\,|\,\HH_n(k[V]) \neq 0\}.
\end{equation}
However, even for
well-behaved noncommutative algebras Hochschild homology
is often rather degenerate. For example,  
the standard quantized coordinate ring 
$\A:=k_q[SL(N)]$ is for generic $q$ Auslander regular and Cohen-Macaulay 
with global and Gelfand-Kirillov dimension 
equal to the classical dimension $N^2 -1$ of $SL(N)$ \cite{aLeSt}, but 
$\HH_n(\A)=0$ for $n \ge N$ \cite{aFeTs}.
In this note we show that this ``dimension drop'' is overcome  
by passing to Hochschild homology $\H_\ast(\A,\M)$ with coefficients in
a suitable bimodule $\M$.

The cosemisimple Hopf algebra 
structure on $\A$ determines the Haar functional 
$h : \A \rightarrow k$ which is left and
right invariant under the coaction of $\A$ on itself, 
and there is a unique automorphism $ \smod \in \mathrm{Aut} (\A)$, the
so-called modular automorphism, such that $h(xy)=h(\sigma (y) x)$ 
for all $x,y \in \A$ (see \cite{bKS}, Section~11.3). 
The crucial coefficient bimodule $\M$ is then  
$\A_\sigma$ which is $\A$ as $k$-vector space with
bimodule structure $x \triangleright y \triangleleft z:=xy \sigma(z)$. 
Our main result is:
\begin{thm}\label{main} 
There is an isomorphism of $k$-vector spaces 
$\H_{N^2-1}(\A,\A_\sigma) \cong k$.
\end{thm}
 For $N=2$ this was shown by explicit calculation in \cite{aKraehmer3}.
The proof for arbitrary $N$ given below relies on the 
following analogue of Poincar\'e duality for 
Hochschild (co)homology proven by van
den Bergh:
\begin{thm} \cite{avdB} \label{vdb}  
Let $\X$ be a smooth 
algebra such that there exists $d_\X \in \mathbb{N} $ 
with $ \H^n(\X,\X^e)=0$ for $n \neq d_\X$, and that 
$\U_\X:=\H^{d_\X}(\X,\X^e)$ is an invertible $\X$-bimodule. 
Then for every $\X$-bimodule $\M$
we have
\begin{equation}
\label{duality}
\H^n(\X,\M) \cong \H_{d_\X-n}(\X,\U_\X \otimes_\X \M).
\end{equation}   
\end{thm}
Here $\X^e:=\X \otimes \X^\mathrm{op}$ is the enveloping algebra of $\X$
(throughout this paper an unadorned $ \otimes $ means 
tensor product over $k$), so 
the Hochschild homology and 
cohomology groups of $\X$ with coefficients in $\M$ are 
$\H_n(\X,\M)=\mathrm{Tor}_n^{\X^e}(\M,\X)$ and 
$\H^n(\X,\M)=\mathrm{Ext}_{\X^e}^n(\X,\M)$, respectively. 
Following \cite{avdB} (erratum) an algebra $\X$ is called smooth 
if it has finite projective dimension 
$\mathrm{pd}_{\X^e}(\X)=\inf\{n \ge 0\,|\,\H^{n+1}(\X,\cdot) = 0\}$ 
as an $\X^e$-module.
As in \cite{CE} we call
$\mathrm{pd}_{\X^e}(\X)$ the dimension of
$\X$ and denote it by $\Hdim(\X)$.
As pointed out by van den Bergh, $\X$ is smooth
if and only if $\X^e$ has finite global dimension. This follows from 
$\gldim(\X) \le \Hdim(\X) \le \gldim(\X^e)$ and
$\Hdim(\X \otimes \Y) \le \Hdim(\X) + \Hdim(\Y)$ 
(\cite{CE}, Propositions 7.4-7.6), which gives
$\gldim(\X^e) \le \Hdim(\X^e) \le 2 \, \Hdim(\X) \le 2 \, \gldim(\X^e)$. 
In the sequel we say that
an algebra has the Poincar\'e duality property if it satisfies the
assumptions of Theorem~\ref{vdb}.

The principal technical result of this note 
consists of remarking successively that Theorem~\ref{vdb} 
applies to the quantized coordinate rings $\B:=k_q[M(N)]$ 
of $N \times N$-matrices, 
$\C:=k_q[GL(N)]$ and
$\A=k_q[SL(N)]$. Theorem~\ref{main} then 
follows from the well-known fact that the center of 
$\A$ consists only of the scalars.

Our main motivation for studying $\H_\ast(\A,\A_\sigma)$
is the so-called twisted cyclic cohomology  
and its link to covariant differential calculi over quantum groups 
both due to Kustermans, Murphy and Tuset \cite{aKMT}.
Twisted cyclic cohomology is defined by a cyclic object in the
sense of Connes \cite{aCo85} depending on an algebra $\X$ and 
an automorphism $ \sigma $. Its  
underlying simplicial homology is 
$\H_\ast(\X,\X_\sigma)$ (at least when $ \sigma $ 
is diagonalizable, see Proposition 2.1 in \cite{aKraehmer3}). 
The volume forms of covariant differential calculi over quantum groups
define twisted cyclic cocycles, with the appearance of the twisting
automorphism forced by the modular properties of the Haar
functional that replaces the traces of Connes' original construction
\cite{aCo85}. In view of Theorem~\ref{vdb} twisted coefficients appear very 
naturally also for purely homological reasons, and
Theorem~\ref{main} and similar results for 
quantum hyperplanes and Podle\'s quantum spheres 
\cite{aHadfield,si} show that the twist determines as in the classical case 
a unique class of top degree in Hochschild homology. 

\section{Proof of Theorem~\ref{main}}
We first consider the quantized coordinate ring
$\B = k_q [ M (N)]$. Recall that this has generators
$u_{ij}$, $1 \leq i, j \leq N$, with relations 
\begin{eqnarray}\label{Mqn}
&& u_{ik} u_{il}=q u_{il} u_{ik}, \quad
	u_{ik} u_{jk}=q u_{jk} u_{ik}, \quad
	u_{ik} u_{jk}=q u_{jk} u_{ik},\nonumber\\ 
&& u_{jk} u_{jl} =q u_{jl} u_{jk}, \quad
	u_{il} u_{jk} = u_{jk} u_{il},\quad
	u_{ik} u_{jl}  - u_{jl} u_{ik}  = (q - q^{-1} ) u_{il} u_{jk}
\end{eqnarray} 
for all $i < j$, $k <l$. Here $q \in k \setminus \{0\}$ 
is a fixed deformation parameter, assumed not to be a root of
unity.
\begin{prop}\label{piotr}
$\B$ has the Poincar\'e duality property with
$d_\B=N^2$ and $\U_\B=\B_\sigma$, with $\sigma$ defined by
\begin{equation}\label{sigma_Dn}
		\s(u_{ij}) := q^{2(N+1-i-j)} u_{ij}
\end{equation}
\end{prop} 
We will use here and later 
the following K\"unneth-type isomorphism of Cartan and Eilenberg:
\begin{thm}\label{CEthm}\cite{CE}, Theorem~XI.3.1. 
Let $k$ be a field, 
$\A_1,\A_2$ be two left Noetherian $k$-algebras and 
$\M_i,\N_i$ be finitely generated left modules over $\A_i$. Then 
\begin{equation}
		\bigoplus_{i+j = n} 
		\Ext^i_{\A_1} ( \M_1  , \N_1  )  \otimes \Ext^j_{\A_2} ( \M_2  , \N_2  ) 
		\cong
		\Ext^n_{\A_1 \otimes \A_2} 
		(\M_1 \otimes \M_2 , \N_1 \otimes \N_2 ).
\end{equation}
\end{thm}
{\bf Proof of Proposition~\ref{piotr}.}
The claim follows from 
Proposition~2 in \cite{avdB}: As mentioned in \cite{manin} 
it follows from a result of Priddy (\cite{priddy}, Theorem~5.3) that 
$\B$ is a graded Koszul algebra.
By definition the Koszul dual $\B^{!}$ has generators
$\hu_{ij}$ with relations orthogonal to (\ref{Mqn}):
\begin{eqnarray}\label{eswird}
&& {\hu_{ij}}^2 = 0 \quad \forall \; i,j,\quad 
	\hu_{ik} \hu_{il} = -q^{-1} \hu_{il} \hu_{ik}, \quad
	\hu_{ik} \hu_{jk} =-q^{-1} \hu_{jk} \hu_{ik}\nonumber\\
&&\hu_{ik} \hu_{jk} =-q^{-1} \hu_{jk} \hu_{ik}, \quad  \hu_{jk} \hu_{jl} =-q^{-1} \hu_{jl} \hu_{jk},\quad
	\hu_{ik} \hu_{jl} = - \hu_{jl} \hu_{ik}, \nonumber\\ 
&&\hu_{il} \hu_{jk} + \hu_{jk} \hu_{il} = (q^{-1}-q) \hu_{ik} \hu_{jl},
\end{eqnarray}  
where $i < j$, $k <l$. These relations  
imply that the monomials $\hu_{i_1j_1} \cdots \hu_{i_nj_n}$, $n=1,\ldots,N^2$, 
$i_1j_1 \prec \ldots \prec i_nj_n$ with respect to lexicographical
ordering, form a $k$-linear basis, and that 
$\B^!$ is Frobenius with Frobenius functional 
$\hat h : \B^! \rightarrow k$ being projection onto the component   
of the longest basis element 
$\hu_{11} \hu_{12} \cdots \hu_{NN-1} \hu_{NN}$ (that is, for each nonzero
$x \in \B^!$ there exists $y \in \B^!$ with $\hat h(xy) \neq 0$).  
The formula for $ \sigma $ 
follows by straightforward computation 
using the relations (\ref{eswird}).

Smoothness of $\B$ follows from some well-known facts about Koszul
algebras (see e.g. the survey \cite{froe}). First, 
$\mathrm{Tor}^{\B^e}_n(k,k) \cong
\mathrm{Ext}^n_{\B^e}(k,k)$, and by 
Theorem~\ref{CEthm} and Koszulity this can be written as 
$\sum_{i+j=n} \B^!_i \otimes \B^!_j$ (note that 
$\B \cong \B^\mathrm{op}$). Thus 
$\mathrm{Tor}^{\B^e}_n(k,k)=0$ for $n>2 N^2$, hence 
$\Hdim(\B) \le 2N^2$ by \cite{b}, Corollary~{8.7.5}.
\qed

It was shown by Farinati 
that the class of algebras having the Poincar\'e duality property  
is closed under localization \cite{aFari}, Theorem~1.5. The quantized coordinate ring 
$\C=k_q[GL(N)]$ of the general linear group is the localization of
$\B=k_q[M(N)]$ at the central quantum determinant 
\cite{manin}
\begin{equation}
		\mathrm{det}_q=\sum_{\pi \in S_N} 
		(-q)^{l(\pi)} u_{1 \pi(1)} \cdots u_{N \pi(N)},
\end{equation}
with $S_N$ the permutation group on $N$ elements and 
$l(\pi)$ the length of a permutation. This is 
$ \sigma $-invariant, so $ \sigma $ passes
to an automorphism of $\C$, still denoted by $ \sigma $
and given by (\ref{sigma_Dn}). Proposition~\ref{piotr} 
now implies:
\begin{cor}\label{22}
 $\C$ has the Poincar\'e duality property with 
$d_\C=d_\B=N^2$ and $\U_\C=\C \otimes_\B \U_\B=\C_\sigma$.
\end{cor}

The algebra $\A=k_q[SL(N)]$ is
the quotient of $\B$ by the relation $\mathrm{det}_q=1$, 
and again by $ \sigma $-invariance of $ \mathrm{det}_q$,
$ \sigma $ descends to an automorphism of $\A$. Following the strategy 
of Levasseur and Stafford \cite{aLeSt} we will use the isomorphism
$\C \cong \A \otimes D$, where $D:=k[t,t^{-1}]$
to deduce Poincar\'e duality for $\A$ from $\B$ via $\C$. 
This enables us to prove finally:
\begin{prop}
The algebra $\A$ has the Poincar\'e duality property with 
$d_\A=N^2-1$ and $\U_\A=\A_\sigma$.
\end{prop}
{\bf Proof.}
We apply Theorem~\ref{CEthm} with 
$\A_1=\N_1=\A^e$, $\M_1=\A$ and 
$\A_2=\N_2=D^e$, $\M_2=D$.
Since $\A=\A^{\mathrm{op}}$ (the antipode of the standard Hopf 
algebra structure gives an isomorphism) we have 
$\A^e \cong k_q[SL(N) \times SL(N)]$, so it is (both left and right) Noetherian
by \cite{bJoseph}, Proposition~{9.2.2} and further $\A$ is smooth by \cite{g}.
It is elementary to show that 
$D$ satisfies Poincar\'e duality with $d_D=1$, $\U_D=D$.
So $\Ext^n_{\A^e} ( \A  , \A^e )  \otimes D \cong \Ext^{n+1}_{\C^e} (\C ,\C^e)$ for each $n \geq 0$.
So by Corollary~\ref{22} we have
$\Ext^{N^2-1}_{\A^e} ( \A  , \A^e )  \otimes D = \C_\sigma $,
which is $\A_\sigma \otimes D$, and all other 
$\Ext^n_{\A^e} ( \A  , \A^e )$ vanish. The result follows.
\qed

Thus there is an isomorphism 
$\H_{N^2-1}(\A,\A_\sigma) \cong \H^0(\A,\A)$.  
The latter is by definition the center of $\A$, and this 
consists only of the scalars (see e.g.~\cite{bJoseph}, 
Theorem~{9.3.20}). This completes the
proof of Theorem~\ref{main}.

\section{Acknowledgements}
It is a pleasure to thank 
Ken Brown for explaining that Theorem~\ref{main} 
can also be proved via the alternative description of Hochschild
(co)homology of Hopf algebras from \cite{aFeTs}, as in  \cite{bz},
and Marco Farinati, Brad Shelton and Paul Smith for helpful comments and
discussion. We are also very grateful to 
Giselle Rowlinson and St\'ephane 
Launois for their  help with the French translation.

\end{document}